\documentclass[11pt,reqno]{article}
\usepackage[T1]{fontenc}
\usepackage{lmodern,microtype}
\usepackage{amsmath,amssymb,mathtools,amsthm,mathrsfs}
\usepackage{tikz}
\usetikzlibrary{arrows.meta}
\usepackage{caption}
\usepackage{booktabs,tabularx,longtable,array}
\usepackage[margin=1.06in]{geometry}
\usepackage{authblk}
\usepackage[hidelinks]{hyperref}
\usepackage{orcidlink}

\theoremstyle{plain}
\newtheorem{theorem}{Theorem}[section]
\newtheorem{proposition}[theorem]{Proposition}
\newtheorem{lemma}[theorem]{Lemma}
\newtheorem{corollary}[theorem]{Corollary}
\theoremstyle{remark}
\newtheorem{remark}[theorem]{Remark}
\usepackage[nameinlink,noabbrev,sort]{cleveref}
\crefname{theorem}{Theorem}{Theorems}
\Crefname{theorem}{Theorem}{Theorems}
\crefname{proposition}{Proposition}{Propositions}
\Crefname{proposition}{Proposition}{Propositions}
\crefname{lemma}{Lemma}{Lemmas}
\Crefname{lemma}{Lemma}{Lemmas}
\crefname{corollary}{Corollary}{Corollaries}
\Crefname{corollary}{Corollary}{Corollaries}
\crefname{remark}{Remark}{Remarks}
\Crefname{remark}{Remark}{Remarks}
\numberwithin{equation}{section}
\newcommand{\norm}[1]{\left\lVert#1\right\rVert}
\newcommand{\abs}[1]{\left|#1\right|}
\newcolumntype{L}[1]{>{\raggedright\arraybackslash}p{#1}}
\newcolumntype{Y}{>{\raggedright\arraybackslash}X}
\DeclareMathOperator{\diver}{div}
\DeclareMathOperator{\supp}{supp}
\hypersetup{
pdftitle={Magnetic norm inflation below the velocity regularity threshold for viscous non-resistive MHD},
pdfauthor={Yinghui Wang},
pdfsubject={Sharp velocity regularity for viscous non-resistive MHD},
pdfkeywords={non-resistive MHD, norm inflation, ill-posedness, unequal Sobolev regularity}}
\title{Magnetic norm inflation below the velocity regularity threshold\\
for viscous non-resistive MHD}
\author{Yinghui Wang\,\orcidlink{0000-0002-7565-5525}}
\affil{\small MOE-LCSM, School of Mathematics and Statistics,
Hunan Normal University, Changsha 410081, China\\
\href{mailto:yhwangmath@hunnu.edu.cn}{yhwangmath@hunnu.edu.cn}}
\date{}
\begin{document}
\maketitle

\begin{abstract}
We study viscous, non-resistive magnetohydrodynamics on $\mathbb R^d$,
$d\in\{2,3\}$. For every fixed magnetic exponent $s>d/2$ and every
velocity exponent $q<s-1$, we prove norm inflation at the origin in
$H^q\times H^s$. Arbitrarily small divergence-free Schwartz data
generate classical solutions whose magnetic $H^s$ norm becomes
arbitrarily large at arbitrarily small positive times, while the
velocity remains arbitrarily small in $H^q$ throughout the observation
interval. Together with endpoint local well-posedness, this identifies
$q=s-1$ as the sharp lower velocity threshold for fixed $s>d/2$.
The construction uses a symmetric pair of high-frequency velocity
packets and a fixed low-frequency magnetic profile. A Fourier component of their
stretching interaction has a fixed sign and gains the power
$K^{s-1-q}$ on a short diffusive time scale. Subtracting an exact
reference flow and using endpoint Besov estimates controls the full
nonlinear remainder. Quantitative higher-order estimates ensure
existence through the observation time and yield the uniform smallness
of the velocity in the target norm.
\end{abstract}


\section{Introduction}
\subsection{The model and the regularity threshold}\label{subsec:model}
We consider the incompressible viscous, non-resistive
magnetohydrodynamics (MHD) system
\begin{equation}\label{eq:MHD}
\left\{
\begin{aligned}
\partial_tu-\Delta u+u\cdot\nabla u+\nabla p&=B\cdot\nabla B,\\
\partial_tB+u\cdot\nabla B&=B\cdot\nabla u,\\
\diver u=\diver B&=0,\\
(u,B)|_{t=0}&=(u_0,B_0),
\end{aligned}
\right.
\qquad (t,x)\in[0,\infty)\times\mathbb R^d,
\end{equation}
where $d\in\{2,3\}$, $u$ is the fluid velocity, $B$ is the magnetic
field, and $p$ is the pressure. The viscosity has been normalized to
one. The magnetic equation has transport and stretching, but no
diffusion. General mathematical background on incompressible MHD may
be found in \cite{SermangeTemam1983}. Here we ask how much rougher the
initial velocity may be than the initial magnetic field while retaining
continuity of the solution map.

For $s>d/2$, Fefferman, McCormick, Robinson and Rodrigo established
local existence and uniqueness first in $H^s\times H^s$
\cite{FeffermanEtAl2014} and then in
$H^{s-1+\varepsilon}\times H^s$, $0<\varepsilon<1$
\cite{FeffermanEtAl2017}. Li proved local well-posedness at the endpoint
\begin{equation}\label{eq:Li-endpoint}
(u_0,B_0)\in H^{s-1}(\mathbb R^d)\times H^s(\mathbb R^d),
\qquad s>d/2,
\end{equation}
using a commutator estimate \cite{Li2020}. Related local theories in
Besov spaces were developed in
\cite{CheminEtAl2016,LiTanYin2017,LuoYeYin2025}.
The one-derivative gap in \eqref{eq:Li-endpoint} suggests a sharp
compatibility condition between the two initial regularities.

\subsection{Main result}\label{subsec:main}
For real $q,s$, let
\begin{equation}\label{eq:phase-space}
\mathcal X_{q,s}
:=H^q_{\mathrm{div}}(\mathbb R^d)
\times H^s_{\mathrm{div}}(\mathbb R^d),
\end{equation}
equipped with the sum norm; the subscript $\mathrm{div}$ denotes
divergence-free vector fields. For Schwartz data, a classical solution
means the smooth Sobolev solution described in \Cref{lem:smooth-flow}.
We denote its maximal smooth lifespan by $T_{\max}(u_0,B_0)$;
we also write $T_{\max}(z_0)$ for $z_0=(u_0,B_0)$.

\begin{theorem}[Magnetic norm inflation with small velocity]
\label[theorem]{thm:main}
Let $d\in\{2,3\}$, $s>d/2$, and $q<s-1$. For every
$\varepsilon>0$ and $M>0$, there exist divergence-free fields
$u_0,B_0\in\mathscr S(\mathbb R^d;\mathbb R^d)$ and
$0<T<\min\{\varepsilon,T_{\max}(u_0,B_0)\}$ such that
\begin{equation}\label{eq:main-conclusion}
\norm{u_0}_{H^q}+\norm{B_0}_{H^s}<\varepsilon,
\qquad
\sup_{0\le t\le T}\norm{u(t)}_{H^q}<\varepsilon,
\qquad
\norm{B(T)}_{H^s}>M.
\end{equation}
In particular, \eqref{eq:MHD} has norm inflation at the origin in
$\mathcal X_{q,s}$, with inflation in the magnetic component.
\end{theorem}

The conclusion rules out continuous extensions of the classical flow
to a common time interval. The following formulation requires agreement
only during each datum's own smooth lifespan.

\begin{corollary}\label[corollary]{cor:no-continuous-solution-map}
Under the assumptions of \Cref{thm:main}, let $T_0>0$ and let
$\mathcal U$ be any neighborhood of zero in $\mathcal X_{q,s}$.
There is no map
\[
\Phi:\mathcal U\longrightarrow C([0,T_0];\mathcal X_{q,s})
\]
that is continuous at zero and, for every divergence-free Schwartz
datum $z_0\in\mathcal U$, satisfies
\begin{equation}\label{eq:agreement-with-classical-solution}
\Phi(z_0)(t)=(u(t),B(t))
\quad\text{for }0\le t<\min\{T_0,T_{\max}(z_0)\},
\end{equation}
where $(u,B)$ is the corresponding classical solution.
\end{corollary}

\begin{remark}[Scope of the threshold]
For each fixed $s>d/2$, \eqref{eq:Li-endpoint} gives local
well-posedness at $q=s-1$, whereas \Cref{thm:main} gives norm inflation
for every $q<s-1$. Sharpness here concerns this lower velocity
threshold with $s$ fixed. In particular, for every $s>d/2$
there is magnetic norm inflation in
$H^{d/2-1}\times H^s$ while the velocity remains
small in $H^{d/2-1}$. 
\end{remark}

\subsection{Relation to earlier instability results}\label{subsec:earlier}
For equal Sobolev regularity $H^s\times H^s$, $s>d/2$,
Li, Yin and Zhu proved continuous dependence and failure of uniform
continuity on bounded sets \cite{LiYinZhu2023}. The closest
norm-inflation result is due to Chen, Nie and Ye
\cite{ChenNieYe2024}: they proved magnetic norm inflation in
$H^{d/2-1}\times H^{d/2}$ and in the corresponding critical Besov
spaces with dyadic summability exponent greater than one. Their
construction identifies the low--high paraproduct in $B\cdot\nabla u$
as the leading magnetic interaction and uses a multiscale accumulation
at the critical magnetic index.

The present result keeps the magnetic index strictly above $d/2$ and
varies the velocity index. Its main contribution is the full range
$q<s-1$ for each fixed $s>d/2$, including
\[
d/2-1<q<s-1.
\]
In this range both indices lie above the pair $(d/2-1,d/2)$, but their
gap exceeds one derivative. The earlier critical result does not imply
this conclusion by Sobolev embedding, since smallness of the magnetic
datum in $H^{d/2}$ does not give smallness in $H^s$. Here one
high-frequency velocity scale and a fixed magnetic profile suffice: the positive
regularity gap $s-1-q$ yields a power gain, replacing critical
multiscale accumulation. The nonlinear
estimates also keep the velocity small in its target norm throughout
the inflation interval.

In three dimensions, Wan, Ye and Yin
\cite[Theorem~1.2 and Remark~1.2]{WanYeYin2026} proved norm inflation
in $H^\sigma\times H^{\sigma+1}$, $0<\sigma<1/2$, for viscous,
non-resistive MHD. They also record magnetic inflation with a
uniformly bounded velocity; the magnetic indices in that result lie
below $3/2$. Ban, Guo and Guo \cite{BanGuoGuo2026} established
magnetic norm inflation at the endpoint
$\dot B^{-1}_{\infty,1}\times\dot B^0_{\infty,1}$.
Instability near a Couette flow with a constant
magnetic field was studied by Knobel and Zillinger
\cite{KnobelZillinger2025} and by Dolce, Knobel and Zillinger
\cite{DolceKnobelZillinger2026}. Those shear-flow problems concern a
nonzero background on a periodic infinite channel and different time
scales or observables. The result here concerns the zero solution on
$\mathbb R^d$ and the compatibility of two Sobolev indices.

\subsection{Outline of the proof}\label{subsec:outline}
\hypertarget{not:scales}{}Write $r=d/2$ and $\alpha=s-1-q>0$. Place frequency $K\gg1$ in the
velocity datum and frequency of order one in the magnetic datum.
Before multiplying by amplitudes $\lambda$ and $\mu$, respectively,
normalize the profiles so that
\[
\norm{u_K^{\rm hi}}_{H^q}\simeq1,
\qquad \norm{u_K^{\rm hi}}_2\simeq K^{-q},
\qquad \norm{b^{\rm lo}}_{H^s}\simeq1.
\]
\hypertarget{not:response}{}The principal magnetic response is
\[
\mathcal M_K(t)
=\int_0^t\mu b^{\rm lo}\cdot\nabla
e^{\varrho\Delta}(\lambda u_K^{\rm hi})\,\mathrm d\varrho.
\]
On $0\le t\le T=\tau K^{-2}$, with $\tau$ small, the heat multiplier
$e^{-t|\xi|^2}$ is bounded below on the velocity packets. The $H^s$ size predicted by the
interaction is therefore
\begin{equation}\label{eq:intro-frequency-count}
\mu\lambda\tau\,
\underbrace{K^s}_{\text{output norm}}\;
\underbrace{K}_{\text{stretching}}\;
\underbrace{K^{-2}}_{\text{time scale}}\;
\underbrace{K^{-q}}_{\text{input size}}
=\mu\lambda\tau K^\alpha.
\end{equation}
The companion term with the derivative on $b^{\rm lo}$ is smaller by
$K^{-1}$. The polarizations and frequency supports in
\Cref{sec_local} select a scalar Fourier component with a fixed sign,
turning \eqref{eq:intro-frequency-count} into a lower bound.

\hypertarget{not:reference}{}To control the exact solution, we subtract the reference solution
$(\bar u,\bar B)$ generated by $(0,\mu b^{\rm lo})$.
This incorporates the full evolution of the low-frequency datum.
\hypertarget{not:endpoint-scale}{}The perturbation is estimated in
$\mathcal B^{r-1}\times\mathcal B^r$, where
$\mathcal B^\sigma=B^\sigma_{2,1}$, and the velocity has size
\[
m:=\lambda K^{r-s+\alpha}.
\]
The main magnetic response has size $\mu m\tau$ in the selected
$\mathcal B^r$ component. The nonlinear remainder is smaller by
$m+\mu+\tau+K^{-1}$, apart from a rapidly decreasing tail of the
reference flow. The time interpolation gain
$\norm{u-\bar u}_{L^1_T\mathcal B^r}\lesssim m\tau K^{-1}$
compensates the extra derivative in transporting the principal
response. Higher-order estimates provide continuation through $T$ and,
by interpolation, smallness of the velocity in $H^q$.

Section~\ref{sec_notation} collects the analytic estimates.
Section~\ref{sec_local} constructs the profiles and proves the
nonlinear expansion. Section~\ref{sec:norm-inflation} chooses the
parameters and concludes the proof.

\section{Notation and analytic tools}\label{sec_notation}

\subsection{Notation}\label{subsec:notation}

The index in \Cref{tab:notation-index} gives the first introduction in the main text and, when
different, the subsequent defining formula. Symbols previewed in the
abstract are included; standard analytic notation is defined here.

In \eqref{eq:MHD}, $t\ge0$ and $x\in\mathbb R^d$, $d\in\{2,3\}$,
are time and position. We use $\partial_t$, $\partial_a=\partial_{x_a}$,
$\nabla$, $\Delta=\sum_{a=1}^d\partial_a^2$, and $\diver$ for time and
spatial differentiation, gradient, Laplacian, and divergence.
A dot denotes the Euclidean inner product or contraction, and
$\otimes$ the tensor product. We use the unitary Fourier transform
\begin{equation}\label{eq:Fourier-Sobolev}
 \begin{aligned}
 \widehat f(\xi)=\mathcal Ff(\xi)
   &=(2\pi)^{-d/2}\int_{\mathbb R^d}e^{-ix\cdot\xi}f(x)\,\mathrm dx,\\
 \norm{f}_{H^\sigma}
   &=\norm{(1+|\xi|^2)^{\sigma/2}\widehat f}_2.
 \end{aligned}
\end{equation}
Here $\xi$ is the Fourier variable, $i^2=-1$, $\mathcal F^{-1}$ is
the inverse transform, and $\norm{\cdot}_2$ is the $L^2(\mathbb R^d)$
norm. Other $L^p$ norms have the usual meaning; vector and matrix
norms use the Euclidean and Frobenius norms. The Sobolev spaces first
occur in Section~\ref{subsec:model}; $H^\sigma_{\rm div}$ in
\eqref{eq:phase-space} is their distributionally divergence-free
subspace, and $\mathscr S$ in \Cref{thm:main} is the Schwartz class.
The notation $C$, $C_c^\infty$, and $L^\rho_{\rm loc}$ denotes,
respectively, continuous functions, smooth compactly supported
functions, and local time integrability, with the displayed domains
and target spaces. These occur in \Cref{cor:no-continuous-solution-map},
Section~\ref{subsec:profiles}, and \eqref{eq:smooth-flow-class}.

For a symbol $a$, set $a(D)f=\mathcal F^{-1}(a\widehat f)$, with
$D=-i\nabla$. The heat semigroup $e^{t\Delta}$, first used in
\hyperlink{not:response}{Section~\ref{subsec:outline}}, has symbol
$e^{-t|\xi|^2}$. The Leray projector and its symbol are
\begin{equation}\label{eq:Leray-symbol}
 \mathbb P=\Pi(D),\qquad
 \Pi(\xi)=I_d-\frac{\xi\otimes\xi}{|\xi|^2}\quad(\xi\ne0),
 \qquad \Pi(0)=I_d,
\end{equation}
where $I_d$ is the identity matrix. Operator commutators are
$[A,B]=AB-BA$; their first use is \eqref{eq:endpoint-commutator-sum}.

Fix a radial inhomogeneous Littlewood--Paley partition
$\chi_0(\xi)+\sum_{j\ge0}\psi(2^{-j}\xi)=1$, with $\chi_0$ supported
in a ball and $\psi$ in an annulus, and define
\begin{equation}\label{eq:dyadic-blocks}
 \begin{gathered}
 \Delta_{-1}=\chi_0(D),\qquad \Delta_j=\psi(2^{-j}D)\quad(j\ge0),
 \qquad S_j=\sum_{-1\le k<j}\Delta_k,\\
 \widetilde\Delta_j=\Delta_{j-1}+\Delta_j+\Delta_{j+1}.
 \end{gathered}
\end{equation}
Blocks below $-1$ are zero; the unindexed $\Delta$ is the Laplacian.
The abbreviation $\mathcal B^\sigma$, introduced in
\hyperlink{not:endpoint-scale}{Section~\ref{subsec:outline}}, and
its time norms are
\begin{equation}\label{eq:Besov-time-norms}
 \begin{aligned}
 \mathcal B^\sigma&=B^\sigma_{2,1}(\mathbb R^d),&
 \norm{f}_{\mathcal B^\sigma}
   &=\sum_{j\ge-1}2^{j\sigma}\norm{\Delta_jf}_2,\\
 &&\norm{f}_{\widetilde L^\rho(I;\mathcal B^\sigma)}
   &=\sum_{j\ge-1}2^{j\sigma}\norm{\Delta_jf}_{L^\rho(I;L^2)}.
 \end{aligned}
\end{equation}
Here $I$ is a time interval and $1\le\rho\le\infty$.
The second line is the Chemin--Lerner norm \cite{CheminLerner1995};
$L^\rho(I;\mathcal B^\sigma)$ is the ordinary Bochner norm, and a
subscript $T$ means $I=(0,T)$. Tonelli's theorem gives
$\widetilde L^1=L^1$. At $t=0$, running $\widetilde L^\infty_t$ norms
take the initial spatial norm; see the continuity argument after
\eqref{eq:CL-interpolation}. The dot in $\dot B^\sigma_{\infty,1}$ in
Section~\ref{subsec:earlier} denotes homogeneous dyadic blocks;
all Besov norms in the proof are inhomogeneous.
We use $\mathcal B^r\hookrightarrow L^\infty$ and
$\mathcal B^\sigma\hookrightarrow H^\sigma$, where
$\hookrightarrow$ denotes continuous embedding and $r=d/2$;
$\mathcal B^r$ is also an algebra.
Bony's paraproduct and remainder \cite{Bony1981,BCD}, used in
\Cref{lem:products}, are
\begin{equation}\label{eq:Bony-decomposition}
 fg=T_fg+T_gf+R(f,g),\qquad
 T_fg=\sum_jS_{j-1}f\,\Delta_jg,\qquad
 R(f,g)=\sum_{|j-k|\le1}\Delta_jf\,\Delta_kg.
\end{equation}

The symbols $\lesssim$ and $\simeq$, first used in
Section~\ref{subsec:outline}, mean inequality up to a constant and
two-sided comparability; $\supp$ denotes support. Generic positive
constants $c,C$ may change between estimates and may depend on fixed
indices, cutoffs, and profiles, but never on $K,\lambda,\mu,\tau$.
Subscripts on $C$ record additional fixed-index dependence.
In the analytic lemmas, $f,g,v,w,F,F_1,F_2$ denote generic fields;
$U,W$ are transport and stretching fields, and subscripts $0,a$
indicate initial values at those times
(\Cref{lem:products,lem:heat,lem:transport,lem:tame-transport}).
The letters $j,k$ are summation indices, $a,b$ are locally specified
indices or interval endpoints, $\sigma,\sigma_0,\sigma_1$ are
regularity orders, and $\varrho$ is a time-integration variable.
These local roles are determined by the displayed ranges.

\begingroup
\small
\setlength{\tabcolsep}{5pt}
\renewcommand{\arraystretch}{1.12}
\begin{longtable}{@{}L{0.24\textwidth}L{0.73\textwidth}@{}}
\caption{Notation index.}\label{tab:notation-index}\\
\toprule
\textbf{Notation} & \textbf{Meaning and first introduction / definition}\\
\midrule
\endfirsthead
\toprule
\textbf{Notation} & \textbf{Meaning and first introduction / definition (continued)}\\
\midrule
\endhead
\bottomrule
\endfoot
$u,B,p;\ u_0,B_0$
& Velocity, magnetic field, pressure, and initial fields:
  \eqref{eq:MHD}.\\
$q,s;\ r,\alpha$
& Velocity and magnetic indices: \eqref{eq:phase-space} and
  Section~\ref{subsec:model}; $r=d/2$, $\alpha=s-1-q$:
  \hyperlink{not:scales}{Section~\ref{subsec:outline}}.
  Section~\ref{sec_local} uses $\alpha>0$;
  Section~\ref{sec:norm-inflation} imposes $\alpha<s-r$.\\
$\mathcal X_{q,s},z_0;\ T_{\max}$
& Divergence-free data space with the sum norm, datum $(u_0,B_0)$,
  and its maximal smooth lifespan: \eqref{eq:phase-space} and
  the paragraph following it; solution class in
  \Cref{lem:smooth-flow}.\\
$\varepsilon,M,T$
& Tolerance, inflation target, and observation time:
  \Cref{thm:main}. In Section~\ref{subsec:model}, $\varepsilon$ is a
  regularity loss; $T=\tau K^{-2}$ first appears in
  \hyperlink{not:response}{Section~\ref{subsec:outline}}.\\
$\mathcal U,\Phi,T_0$
& Neighborhood of zero, candidate continuous flow map, and its
  common time interval: \Cref{cor:no-continuous-solution-map}.\\
$K,\lambda,\mu,\tau,m$
& Dyadic frequency, velocity/magnetic amplitudes, time fraction,
  and derived size $m=\lambda K^{r-s+\alpha}$:
  \hyperlink{not:scales}{Section~\ref{subsec:outline}},
  \hyperlink{not:response}{response} ($\tau$), and
  \hyperlink{not:endpoint-scale}{endpoint estimate} ($m$);
  data in \eqref{eq:data-and-parameters}.\\
$u_K^{\rm hi},b^{\rm lo}$
& High-frequency velocity and fixed low-frequency magnetic profiles:
  \hyperlink{not:scales}{Section~\ref{subsec:outline}};
  formulas \eqref{eq:u-high-profile}--\eqref{eq:b-low-profile}.\\
$e_1,e_2,\mathbb B(\xi_0,R)$
& Coordinate unit vectors and the Fourier ball with center $\xi_0$
  and radius $R$: opening of Section~\ref{sec_local}.\\
$\kappa,\varphi,\chi_{\rm loc}$
& Fixed packet radius, even profile bump, and output cutoff:
  Section~\ref{subsec:profiles}, before
  \eqref{eq:packet-polarizations}, \eqref{eq:u-high-profile}, and
  \eqref{eq:Qloc-def}, respectively.\\
$\theta,\zeta$
& Low- and high-frequency convolution variables in
  \eqref{eq:packet-polarizations}; $\zeta=\xi-\theta$ in the
  proof of \Cref{lem:M-lower}.\\
$\mathcal Q_K$
& Annular Fourier localization: \Cref{lem:localization};
  the selected output packet is specified by \eqref{eq:Qloc-def}.\\
$\mathcal M_K$
& Leading magnetic stretching response:
  \hyperlink{not:response}{Section~\ref{subsec:outline}};
  defining integral \eqref{eq:M-def}.\\
$(\bar u,\bar B),\ (v,h)$
& Exact reference solution with datum $(0,\mu b^{\rm lo})$:
  \hyperlink{not:reference}{Section~\ref{subsec:outline}} and
  \Cref{lem:background}; perturbations
  $v=u-\bar u$, $h=B-\bar B$: \eqref{eq:perturbation-definition}.\\
$v_{\rm lin},v_{\rm rem}$
& Heat evolution of the velocity datum and remaining velocity
  perturbation: \eqref{eq:velocity-splitting}.\\
$\mathcal R_K$
& Magnetic remainder $h-\mathcal M_K$:
  \eqref{eq:magnetic-remainder-definition}; distinct from the
  bilinear Bony remainder $R(f,g)$ in \eqref{eq:Bony-decomposition}.\\
$\ell,k;\ N,J$
& Integer continuation order and higher Sobolev order:
  \Cref{lem:smooth-flow} and its proof. The tame order $N$ starts in
  \Cref{lem:products} and is fixed as an integer in
  \Cref{prop:critical}; $J$ is the auxiliary reference order in
  the proof of \Cref{lem:background}.\\
$L,\gamma,\delta$
& Tail order: \Cref{lem:localization}; interpolation weight
  $0\le\gamma\le1$: \eqref{eq:CL-interpolation}; final decay exponent:
  \eqref{eq:delta-choice}.\\
$\eta_N,\eta_{N,L}$
& Smallness thresholds, fixed independently of $K$ after the indicated
  orders: \Cref{prop:critical,prop:expansion}.\\
$\mathcal E_{u,J},\mathcal E_{B,J}$
& Reference velocity and magnetic running norms:
  \eqref{eq:reference-running-norms}.\\
$X,D,Y,E$
& Endpoint velocity supremum, velocity dissipation integral, magnetic
  perturbation supremum, and velocity-remainder norm:
  \eqref{eq:endpoint-running-norms}. Here $D(t)$ is a scalar norm,
  whereas $D$ in $a(D)$ is the differential operator.\\
$\mathcal V_N,\mathcal H_N$
& Higher-order velocity and magnetic perturbation controls:
  \eqref{eq:high-running-norms}.\\
$A_u,A_B;\ A_v,A_h$
& Fixed reference bootstrap coefficients:
  \hyperlink{not:reference-barriers}{proof of \Cref{lem:background}};
  fixed perturbation coefficients: \eqref{eq:perturbation-barriers}.\\
$c_0,C_0,C_{\rm tr},C_*$
& Fixed response lower-bound, linear-packet, transport, and derived
  bootstrap constants, respectively:
  \eqref{eq:M-lower}, \eqref{eq:linear-truncated},
  \eqref{eq:magnetic-perturbation-estimate},
  \eqref{eq:remainder-bootstrap}.\\
$n,z_{0,n},B_n,T_n$
& Sequence index, initial data, associated magnetic solutions, and
  observation times:
  \hyperlink{not:corollary-sequence}{proof of
  \Cref{cor:no-continuous-solution-map}}.\\
\end{longtable}
\endgroup

\subsection{Products and evolution estimates}

\begin{lemma}[Endpoint and tame products]\label[lemma]{lem:products}
The following estimates hold:
\begin{align}
 \norm{v\cdot\nabla w}_{\mathcal B^{r-1}}
 &\lesssim \norm{v}_{\mathcal B^{r-1}}
             \norm{w}_{\mathcal B^{r+1}},
             \label{eq:NS-product}\\
 \norm{F_1\cdot\nabla F_2}_{\mathcal B^{r-1}}
 &\lesssim \norm{F_1}_{\mathcal B^r}\norm{F_2}_{\mathcal B^r},
 \qquad \diver F_1=0,
             \label{eq:mag-product}\\
 \norm{F\cdot\nabla v}_{\mathcal B^r}
 &\lesssim \norm{F}_{\mathcal B^r}\norm{v}_{\mathcal B^{r+1}}.
             \label{eq:stretch-product}
\end{align}
For $N>r+1$, one also has
\begin{align*}
 \norm{v\cdot\nabla w}_{\mathcal B^{N-1}}
 &\le C_N\bigl(\norm{v}_{\mathcal B^{r-1}}
                    \norm{w}_{\mathcal B^{N+1}}
             +\norm{v}_{\mathcal B^{N-1}}
                    \norm{w}_{\mathcal B^{r+1}}\bigr),\\
 \norm{F\cdot\nabla v}_{\mathcal B^N}
 &\le C_N\bigl(\norm{F}_{\mathcal B^r}\norm{v}_{\mathcal B^{N+1}}
             +\norm{F}_{\mathcal B^N}\norm{v}_{\mathcal B^{r+1}}\bigr),\\
 \norm{F_1\cdot\nabla F_2}_{\mathcal B^{N-1}}
 &\le C_N\bigl(\norm{F_1}_{\mathcal B^r}\norm{F_2}_{\mathcal B^N}
             +\norm{F_1}_{\mathcal B^N}\norm{F_2}_{\mathcal B^r}\bigr),
 \quad \diver F_1=0.
\end{align*}
Finally, for $N>r$,
\begin{equation}\label{eq:high-algebra-product}
 \norm{f\cdot\nabla g}_{\mathcal B^N}
 \le C_N\norm{f}_{\mathcal B^N}\norm{g}_{\mathcal B^{N+1}}.
\end{equation}
\end{lemma}

\begin{proof}
Use Bony's decomposition \eqref{eq:Bony-decomposition}. The resulting
product law is
\begin{equation}\label{eq:critical-product-law}
 \norm{fg}_{\mathcal B^{a+b-r}}
 \lesssim\norm{f}_{\mathcal B^a}\norm{g}_{\mathcal B^b},
 \qquad a,b\le r,\quad a+b>0.
\end{equation}
At the lowest endpoint used here, $d=2$ and $(a,b)=(0,1)$,
the remainder is summable. With the enlarged blocks from
\eqref{eq:dyadic-blocks}, Bernstein's inequality gives
\[
 \sum_j\norm{\Delta_jR(f,g)}_2
 \lesssim\sum_j\sum_{k\ge j-c}2^{j-k}
      \norm{\Delta_kf}_2\,
      2^k\norm{\widetilde\Delta_kg}_2
 \lesssim\norm{f}_{\mathcal B^0}\norm{g}_{\mathcal B^1}.
\]
Here $c$ is a fixed integer determined by the partition. The
paraproducts obey the same bound. Thus \eqref{eq:NS-product} follows
from \eqref{eq:critical-product-law} with $(a,b)=(r-1,r)$ in both
dimensions. For \eqref{eq:mag-product}, write
$F_1\cdot\nabla F_2=\sum_{a=1}^d\partial_a(F_1^aF_2)$ and use the
algebra property of $\mathcal B^r$; this property also proves
\eqref{eq:stretch-product}. In the first tame estimate,
$\norm{S_{j-1}v}_\infty\lesssim2^j\norm{v}_{\mathcal B^{r-1}}$
controls the low--high paraproduct, while the other terms place the
high norm on $v$. The remaining tame bounds follow from the usual
algebra estimate with one factor in $\mathcal B^r$, using the same
divergence identity when required. The algebra property of
$\mathcal B^N$ proves \eqref{eq:high-algebra-product}.
\end{proof}

\begin{lemma}[Heat estimate]\label[lemma]{lem:heat}
If $\partial_tw-\Delta w=f$, $w(0)=w_0$, then for $\sigma\in\mathbb R$
and $0<T\le1$,
\begin{equation}\label{eq:heat-estimate}
 \norm{w}_{\widetilde L^\infty_T\mathcal B^\sigma}
 +\norm{w}_{L^1_T\mathcal B^{\sigma+2}}
 \le C_\sigma\bigl(\norm{w_0}_{\mathcal B^\sigma}
                    +\norm{f}_{L^1_T\mathcal B^\sigma}\bigr).
\end{equation}
The same estimate holds with forcing $\mathbb Pf$ and on translated
intervals of length at most one.
\end{lemma}

\begin{proof}
For $j\ge0$, the dyadic energy inequality is
$\partial_t\norm{\Delta_jw}_2+c2^{2j}\norm{\Delta_jw}_2
\le\norm{\Delta_jf}_2$. Integrate it, take the time supremum, and
sum with weight $2^{j\sigma}$. The low block is bounded without
coercivity, using $\norm{\Delta_{-1}w}_{L^1_TL^2}
\le T\norm{\Delta_{-1}w}_{L^\infty_TL^2}$.
The projector $\mathbb P$ is bounded on every space in
\eqref{eq:heat-estimate} by its bounded $L^2$ multiplier symbol.
\end{proof}

\begin{lemma}[Endpoint transport estimate]\label[lemma]{lem:transport}
Let $\diver U=0$ and suppose
\begin{equation}\label{eq:abstract-transport}
 \partial_tF+U\cdot\nabla F=F\cdot\nabla W+g,
 \qquad F(0)=F_0.
\end{equation}
Then
\begin{align}
 \norm{F}_{\widetilde L^\infty_T\mathcal B^r}
 &\le C\exp\!\left(C\int_0^T
       (\norm{U}_{\mathcal B^{r+1}}+\norm{W}_{\mathcal B^{r+1}})
       \,\mathrm dt\right)\notag\\
 &\qquad\times\bigl(\norm{F_0}_{\mathcal B^r}
                    +\norm{g}_{L^1_T\mathcal B^r}\bigr).
 \label{eq:transport-estimate}
\end{align}
No divergence condition on $F$ or $W$ is needed.
\end{lemma}

\begin{proof}
The dyadic commutator bound
\begin{equation}\label{eq:endpoint-commutator-sum}
 \sum_{j\ge-1}2^{jr}\norm{[U\cdot\nabla,\Delta_j]F}_2
 \lesssim\norm{U}_{\mathcal B^{r+1}}\norm{F}_{\mathcal B^r}
\end{equation}
follows from the paraproduct decomposition; see
\cite[Lemma~2.100 and Theorem~3.14]{BCD}. Apply the dyadic $L^2$
energy estimate to \eqref{eq:abstract-transport}, use
\eqref{eq:stretch-product}, and take the time supremum in each block
before summing. Gronwall's inequality yields
\eqref{eq:transport-estimate} in the stated Chemin--Lerner norm.
\end{proof}

\begin{lemma}[Tame transport estimate]\label[lemma]{lem:tame-transport}
Let $N>r+1$, $\diver U=0$, and
$\partial_tF+U\cdot\nabla F=F\cdot\nabla W+g$ on $(a,b)$, with
$F(a)=F_a$. For $a<t\le b$,
\begin{align}
 \norm{F}_{\widetilde L^\infty(a,t;\mathcal B^N)}
 &\le C_N\exp\!\left(C_N\int_a^t
   (\norm{U}_{\mathcal B^{r+1}}+\norm{W}_{\mathcal B^{r+1}})
   \,\mathrm d\varrho\right)\notag\\
 &\quad\times\left[\norm{F_a}_{\mathcal B^N}
    +\norm{g}_{L^1(a,t;\mathcal B^N)}\right.\notag\\
 &\hspace{7em}\left.
    +\int_a^t\norm{F}_{\mathcal B^r}
       (\norm{U}_{\mathcal B^{N+1}}+\norm{W}_{\mathcal B^{N+1}})
       \,\mathrm d\varrho\right].
 \label{eq:tame-transport-estimate}
\end{align}
\end{lemma}

\begin{proof}
The tame form of \eqref{eq:endpoint-commutator-sum} is
\[
 \sum_j2^{jN}\norm{[U\cdot\nabla,\Delta_j]F}_2
 \le C_N\bigl(\norm{U}_{\mathcal B^{r+1}}\norm{F}_{\mathcal B^N}
             +\norm{U}_{\mathcal B^{N+1}}\norm{F}_{\mathcal B^r}\bigr).
\]
It follows from the same decomposition, placing the high norm on
either factor. Combine it with the tame stretching bound in
\Cref{lem:products} and apply the preceding dyadic argument.
\end{proof}

\begin{lemma}[Tame parabolic estimate]\label[lemma]{lem:tame-parabolic}
Let $N>r+1$ and $I=(a,b)$ with $b-a\le1$. If
$\partial_tw-\Delta w+U\cdot\nabla w=f$ and $w(a)=w_a$, then
\begin{align}
 &\norm{w}_{\widetilde L^\infty(I;\mathcal B^{N-1})}
  +\norm{w}_{L^1(I;\mathcal B^{N+1})}\notag\\
 &\quad\le C_N\left[\norm{w_a}_{\mathcal B^{N-1}}
   +\norm{f}_{L^1(I;\mathcal B^{N-1})}\right.\notag\\
 &\hspace{6em}
   +\norm{U}_{\widetilde L^\infty(I;\mathcal B^{r-1})}
       \norm{w}_{L^1(I;\mathcal B^{N+1})}\notag\\
 &\hspace{6em}\left.
   +\int_I\norm{U}_{\mathcal B^{N-1}}\norm{w}_{\mathcal B^{r+1}}
       \,\mathrm dt\right].
 \label{eq:tame-parabolic-estimate}
\end{align}
This remains valid for $\mathbb P(U\cdot\nabla w)$ in place of
$U\cdot\nabla w$.
\end{lemma}

\begin{proof}
Apply \Cref{lem:heat} to $f-U\cdot\nabla w$ and use the first tame
product estimate in \Cref{lem:products}.
\end{proof}

\subsection{Smooth solutions, interpolation, and localization}

\begin{lemma}[Classical solutions and continuation]\label[lemma]{lem:smooth-flow}
For an integer $\ell>r+1$ and divergence-free data
$(u_0,B_0)\in H^\ell\times H^\ell$, system \eqref{eq:MHD} has a unique
maximal solution on $[0,T_{\max})$ with
\begin{equation}\label{eq:smooth-flow-class}
 u\in C([0,T_{\max});H^\ell)\cap
       L^2_{\rm loc}([0,T_{\max});H^{\ell+1}),
 \qquad B\in C([0,T_{\max});H^\ell).
\end{equation}
If $T_{\max}<\infty$, a bound
\begin{equation}\label{eq:smooth-continuation-criterion}
 \sup_{0\le t<T_{\max}}
    (\norm{u(t)}_{H^\ell}+\norm{B(t)}_{H^\ell})<\infty
\end{equation}
allows continuation past $T_{\max}$. For Schwartz initial data, all
Sobolev orders persist on every compact subinterval of this lifespan.
\end{lemma}

\begin{proof}
The Sobolev local theory gives \eqref{eq:smooth-flow-class} with a
lifespan bounded below solely in terms of $\ell$ and the data norms
\cite[Theorem~1.1 and Section~3]{FeffermanEtAl2014}. Restarting at
times approaching $T_{\max}$ proves the continuation assertion.
For an integer $k\ge\ell$, the divergence-free transport and magnetic
cancellations give
\[
 \frac{\mathrm d}{\mathrm dt}
   (\norm{u}_{H^k}^2+\norm{B}_{H^k}^2)
 +\norm{\nabla u}_{H^k}^2
 \le C_k(\norm{\nabla u}_\infty+\norm{\nabla B}_\infty)
       (\norm{u}_{H^k}^2+\norm{B}_{H^k}^2).
\]
The coefficient is bounded on compact subintervals by the $H^\ell$
norms. Approximation and Gronwall's inequality at each $k$ give the
final statement.
\end{proof}

Dyadic H\"older and time Cauchy--Schwarz inequalities give
\begin{equation}\label{eq:time-interpolation}
 \norm{v}_{L^1_T\mathcal B^r}
 \le T^{1/2}
    \norm{v}_{\widetilde L^\infty_T\mathcal B^{r-1}}^{1/2}
    \norm{v}_{L^1_T\mathcal B^{r+1}}^{1/2}.
\end{equation}
We also use, for $0\le\gamma\le1$,
\begin{equation}\label{eq:CL-interpolation}
 \norm{v}_{\widetilde L^\infty_T
              \mathcal B^{(1-\gamma)\sigma_0+\gamma\sigma_1}}
 \le\norm{v}_{\widetilde L^\infty_T\mathcal B^{\sigma_0}}^{1-\gamma}
      \norm{v}_{\widetilde L^\infty_T\mathcal B^{\sigma_1}}^\gamma.
\end{equation}
For a solution with all Sobolev orders, the running norms in these
estimates and in \eqref{eq:heat-estimate} are continuous on every
compact time interval: finite dyadic sums are continuous and a higher
Sobolev bound makes their tails uniformly small. With the initial-value
convention in Section~\ref{subsec:notation}, this also applies at zero
in the stopping-time argument.

\begin{lemma}[Frequency localization]\label[lemma]{lem:localization}
Let $K\ge1$ and let $\mathcal Q_K$ be a Fourier multiplier with
uniformly bounded symbol supported in $cK\le|\xi|\le CK$, for fixed
$0<c<C$. For $\sigma\in\mathbb R$ and $L>0$,
\begin{align}
 \norm{\mathcal Q_Kf}_2
 &\le C_\sigma K^{-\sigma}\norm{f}_{H^\sigma},
 \label{eq:loc-Hs}\\
 K^\sigma\norm{\mathcal Q_Kf}_2
 &\le C_\sigma\norm{f}_{\mathcal B^\sigma},
 \label{eq:loc-Bs}\\
 K^\sigma\norm{\mathcal Q_Kf}_2
 &\le C_{\sigma,L}K^{-L}\norm{f}_{\mathcal B^{\sigma+L}}.
 \label{eq:loc-Bs-tail}
\end{align}
\end{lemma}

\begin{proof}
Only a uniformly bounded number of shells with $2^j\simeq K$
intersect the multiplier support. Apply Plancherel's theorem and the
definitions of the Sobolev and Besov norms.
\end{proof}

\section{Construction and localized expansion}\label{sec_local}

Fix $s>r=d/2$ and $\alpha>0$. Throughout this section,
$K\ge16$ is dyadic, $e_1,e_2$ are coordinate vectors, and
$\mathbb B(\xi_0,R)$ denotes a Fourier ball.

\subsection{Frequency-localized initial data}\label{subsec:profiles}

Choose $0<\kappa<10^{-3}$ so small that, uniformly for $K\ge16$,
$\theta\in\mathbb B(e_2,\kappa)$, and
$\zeta\in\mathbb B(Ke_1,\kappa)$,
\begin{equation}\label{eq:packet-polarizations}
 \frac K2\le|\zeta|\le2K,\qquad
 (\Pi(\theta)e_1)\cdot\zeta\ge\frac K2,\qquad
 e_2\cdot\Pi(\zeta)e_2\ge\frac12.
\end{equation}
These bounds follow from continuity of the projector symbol near
the two central directions. Let
$\varphi\in C_c^\infty(\mathbb B(0,\kappa))$ be even and
nonnegative, with $\varphi\ge1$ on $\mathbb B(0,\kappa/2)$.
Define
\begin{align}
 \widehat{u_K^{\mathrm{hi}}}(\xi)
 &=K^{-(s-1-\alpha)}
   [\varphi(\xi-Ke_1)+\varphi(\xi+Ke_1)]\Pi(\xi)e_2,
 \label{eq:u-high-profile}\\
 \widehat{b^{\mathrm{lo}}}(\xi)
 &=[\varphi(\xi-e_2)+\varphi(\xi+e_2)]\Pi(\xi)e_1.
 \label{eq:b-low-profile}
\end{align}
Both fields are real, divergence-free, and Schwartz: their Fourier
transforms are smooth, compactly supported away from zero, and real
and even. The magnetic polarization is parallel to the central
velocity frequency, so stretching contributes a factor of order $K$.

\begin{figure}[htbp]
\centering
\captionsetup{font=footnotesize,labelfont=bf}
\begin{tikzpicture}[x=0.75cm,y=0.82cm,>=Latex,
  every node/.style={font=\small}]
  \begin{scope}[xshift=-3.85cm]
    \draw[->] (-4.1,0)--(4.25,0) node[right] {$\xi_1$};
    \draw[->] (0,-2.05)--(0,2.15) node[above] {$\xi_2$};
    \fill (0,0) circle (1.2pt);
    \filldraw[fill=gray!12,draw=black,thick] (-3.15,0) circle (0.29);
    \filldraw[fill=gray!12,draw=black,thick] ( 3.15,0) circle (0.29);
    \filldraw[fill=white,draw=black,densely dashed]
      (0,-1.30) circle (0.29);
    \filldraw[fill=white,draw=black,densely dashed]
      (0, 1.30) circle (0.29);
    \foreach \x/\y in {-3.15/0,3.15/0,0/-1.30,0/1.30}
      {\fill (\x,\y) circle (1.25pt);}
    \node[below] at (-3.15,-0.33) {$(-K,0)$};
    \node[below] at ( 3.15,-0.33) {$(K,0)$};
    \node[left] at (-0.32, 1.30) {$(0,1)$};
    \node[left] at (-0.32,-1.30) {$(0,-1)$};
    \node[align=center,font=\footnotesize] at (0,-2.42)
      {Input packets: velocity (solid), magnetic field (dashed)};
  \end{scope}
  \begin{scope}[xshift=3.85cm]
    \draw[->] (-4.1,0)--(4.25,0) node[right] {$\xi_1$};
    \draw[->] (0,-2.05)--(0,2.15) node[above] {$\xi_2$};
    \fill (0,0) circle (1.2pt);
    \foreach \x/\y in {-3.15/-1.30,-3.15/1.30,3.15/-1.30}
      {\filldraw[fill=gray!18,draw=black!55] (\x,\y) circle (0.29);
       \fill (\x,\y) circle (1.25pt);}
    \filldraw[fill=gray!35,draw=black,very thick]
      (3.15,1.30) circle (0.29);
    \fill (3.15,1.30) circle (1.25pt);
    \draw[black,dashed,thick] (3.15,1.30) circle (0.50);
    \node[below left] at (-3.18,-1.63) {$(-K,-1)$};
    \node[above left] at (-3.18, 1.63) {$(-K,1)$};
    \node[below right] at (3.18,-1.63) {$(K,-1)$};
    \node[above right] at (3.18, 1.63) {$(K,1)$};
    \node[fill=white,inner sep=1pt,align=center]
      at (1.62,0.12) {$\mathcal Q_K$ selects\\$(K,1)$};
    \draw[->]
      (2.18,0.38) .. controls (2.48,0.56) and (2.70,0.82) .. (2.91,1.03);
    \node[align=center,font=\footnotesize] at (0,-2.42)
      {Convolution outputs and the selected packet};
  \end{scope}
\end{tikzpicture}
\caption{Fourier support geometry, drawn schematically and not to
scale. In dimension three the disks represent projections of the
corresponding balls. All packet radii are fixed independently of $K$.}
\label{fig:support-geometry}
\end{figure}
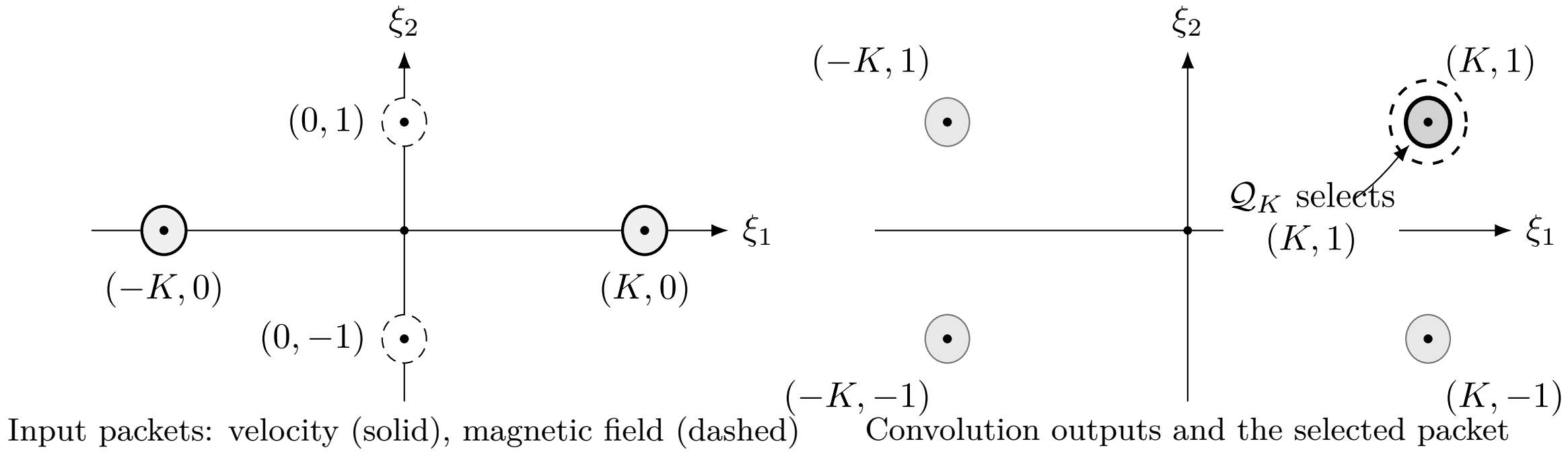

Choose $\chi_{\mathrm{loc}}\in C_c^\infty(\mathbb R^d)$ with
$0\le\chi_{\mathrm{loc}}\le1$, supported in
$\mathbb B(0,4\kappa)$ and equal to one on
$\mathbb B(0,3\kappa)$, and set
\begin{equation}\label{eq:Qloc-def}
 \widehat{\mathcal Q_Kf}(\xi)
 :=\chi_{\mathrm{loc}}(\xi-Ke_1-e_2)\widehat f(\xi).
\end{equation}
For every $t\ge0$, the convolution support satisfies
\begin{equation}\label{eq:four-output-packets}
 \supp\mathcal F\bigl(b^{\mathrm{lo}}\cdot\nabla
                   e^{t\Delta}u_K^{\mathrm{hi}}\bigr)
 \subset
 \bigcup_{a,b\in\{-1,1\}}
 \mathbb B(aKe_1+be_2,2\kappa).
\end{equation}
The four balls are disjoint. The multiplier $\mathcal Q_K$ is one
on the ball centered at $Ke_1+e_2$ and vanishes on the other three.
Thus it selects precisely the interaction between the packets near
$e_2$ and $Ke_1$; see \Cref{fig:support-geometry}. Its symbol is
uniformly bounded and supported where $|\xi|\simeq K$, so
\Cref{lem:localization} applies uniformly in $K$.

\begin{lemma}[Profile norms]\label[lemma]{lem:profile-norms}
For every fixed $\sigma\in\mathbb R$,
\begin{align}
 \norm{u_K^{\mathrm{hi}}}_{H^\sigma}
 &\simeq\norm{u_K^{\mathrm{hi}}}_{\mathcal B^\sigma}
 \simeq K^{\sigma-s+1+\alpha},\notag\\
 \norm{b^{\mathrm{lo}}}_{H^\sigma}
 &\simeq\norm{b^{\mathrm{lo}}}_{\mathcal B^\sigma}
 \simeq1,
 \label{eq:profile-norms}
\end{align}
with constants independent of $K$.
\end{lemma}

\begin{proof}
The velocity packets have fixed volume, lie at $|\xi|\simeq K$,
and have polarization bounded below on the smaller balls where
$\varphi\ge1$. Plancherel's theorem gives
$\norm{u_K^{\mathrm{hi}}}_2\simeq K^{-s+1+\alpha}$.
Only a bounded number of dyadic shells meet their support, and all
have weight comparable to $K^\sigma$. The same argument on the
fixed balls around $\pm e_2$ proves the magnetic estimates.
\end{proof}

For positive amplitudes $\lambda,\mu$ and $0<\tau\le1/4$, set
\begin{equation}\label{eq:data-and-parameters}
 u_0=\lambda u_K^{\mathrm{hi}},\qquad
 B_0=\mu b^{\mathrm{lo}},\qquad
 T=\tau K^{-2},\qquad
 m=\lambda K^{r-s+\alpha}.
\end{equation}
Here $m$ is the derived endpoint velocity size,
$\norm{u_0}_{\mathcal B^{r-1}}\simeq m$, and $\tau$ is the
fraction of the diffusive time scale. Define the leading magnetic
response by
\begin{equation}\label{eq:M-def}
 \mathcal M_K(t)
 :=\int_0^t\mu b^{\mathrm{lo}}\cdot\nabla
             e^{\varrho\Delta}(\lambda u_K^{\mathrm{hi}})
             \,\mathrm d\varrho.
\end{equation}
The companion term with the derivative on $b^{\mathrm{lo}}$ will be
included in the remainder.

\begin{lemma}[Selected magnetic response]\label[lemma]{lem:M-lower}
There are $c_0,C>0$, independent of $K,\lambda,\mu,\tau$, such that
\begin{align}
 K^r\norm{\mathcal Q_K\mathcal M_K(T)}_2
 &\ge c_0\mu m\tau,
 \label{eq:M-lower}\\
 \norm{\mathcal M_K}_{\widetilde L^\infty_T\mathcal B^r}
 +K^{-1}\norm{\mathcal M_K}_{\widetilde L^\infty_T\mathcal B^{r+1}}
 &\le C\mu m\tau.
 \label{eq:M-upper}
\end{align}
\end{lemma}

\begin{proof}
Write $\zeta=\xi-\theta$. By \eqref{eq:four-output-packets}, the
localized convolution involves only
$\theta\in\mathbb B(e_2,\kappa)$ and
$\zeta\in\mathbb B(Ke_1,\kappa)$. With the Fourier normalization fixed above,
\begin{align*}
 \widehat{\mathcal Q_K\mathcal M_K(T)}(\xi)
 &=i(2\pi)^{-d/2}\lambda\mu K^{-s+1+\alpha}
       \chi_{\mathrm{loc}}(\xi-Ke_1-e_2)\\
 &\quad\times\int
       \varphi(\theta-e_2)
       \bigl(\Pi(\theta)e_1\cdot\zeta\bigr)
       \varphi(\zeta-Ke_1)\Pi(\zeta)e_2
       \left(\int_0^Te^{-\varrho|\zeta|^2}\,\mathrm d\varrho\right)
       \,\mathrm d\theta.
\end{align*}
Take the scalar component
$-i\,e_2\cdot\widehat{\mathcal Q_K\mathcal M_K(T)}(\xi)$.
Every factor in its integral is nonnegative by
\eqref{eq:packet-polarizations}. If
$\xi\in\mathbb B(Ke_1+e_2,\kappa/8)$, restrict the integral to
$\theta\in\mathbb B(e_2,\kappa/8)$.
Then $\zeta\in\mathbb B(Ke_1,\kappa/4)$, both profile factors are
at least one, and $\chi_{\mathrm{loc}}=1$. Moreover,
$T|\zeta|^2\le4\tau\le1$, so
\[
 \int_0^Te^{-\varrho|\zeta|^2}\,\mathrm d\varrho\ge e^{-1}T.
\]
The integration region and output ball have fixed positive volume.
Consequently,
\[
 \abs{\widehat{\mathcal Q_K\mathcal M_K(T)}(\xi)}
 \ge c\lambda\mu TK^{2-s+\alpha}
 \quad\text{on }\mathbb B(Ke_1+e_2,\kappa/8).
\]
Plancherel's theorem and \eqref{eq:data-and-parameters} give
\eqref{eq:M-lower}.

For the upper bound, the same convolution formula on all four
packets gives
$\norm{\mathcal M_K(t)}_2\le C\lambda\mu tK^{2-s+\alpha}$.
The support remains in a bounded number of shells at frequency $K$.
Taking the supremum for $t\le T$ and multiplying by $K^r$ or
$K^{r+1}$ yields \eqref{eq:M-upper}.
\end{proof}

\subsection{Reference solution}

To incorporate the self-interaction of the low-frequency magnetic
datum, we subtract its exact evolution under \eqref{eq:MHD}.

\begin{lemma}[Bounds for the reference solution]\label[lemma]{lem:background}
Let $(\bar u,\bar B)$ solve \eqref{eq:MHD} with initial data
$(0,\mu b^{\mathrm{lo}})$. For every fixed $\sigma\ge0$, if $\mu>0$ is
sufficiently small depending on $\sigma$, this solution exists on
$[0,1]$ and, for $0<T\le1$,
\begin{align}
 \norm{\bar u}_{\widetilde L^\infty_T\mathcal B^\sigma}
 +\norm{\bar u}_{L^1_T\mathcal B^{\sigma+2}}
 &\le C_\sigma\mu^2T,\label{eq:bg-u}\\
 \norm{\bar B}_{\widetilde L^\infty_T\mathcal B^\sigma}
 &\le C_\sigma\mu,\label{eq:bg-B-size}\\
 \norm{\bar B-\mu b^{\mathrm{lo}}}_{\widetilde L^\infty_T\mathcal B^\sigma}
 &\le C_\sigma\mu^3T.\label{eq:bg-B-increment}
\end{align}
The constants depend on $\sigma$ and the fixed profile, but not on $K$.
\end{lemma}

\begin{proof}
Choose an integer $J>\max\{r+3,\sigma+2\}$ and consider the maximal
solution in $H^J\times H^J$ supplied by \Cref{lem:smooth-flow},
and denote its lifespan by $T_{\max}$. For
$0<t<\min\{1,T_{\max}\}$, set
\begin{equation}\label{eq:reference-running-norms}
 \mathcal E_{u,J}(t)
 :=\norm{\bar u}_{\widetilde L^\infty_t\mathcal B^J}
   +\norm{\bar u}_{L^1_t\mathcal B^{J+2}},
 \qquad
 \mathcal E_{B,J}(t)
 :=\norm{\bar B}_{\widetilde L^\infty_t\mathcal B^{J+1}}.
\end{equation}
The heat estimate, the high-order algebra bound, and interpolation give
\begin{equation}\label{eq:bg-high-u}
 \mathcal E_{u,J}(t)
 \le C_J\bigl(t^{1/2}\mathcal E_{u,J}(t)^2
                  +t\mathcal E_{B,J}(t)^2\bigr).
\end{equation}
Indeed, the convection term is bounded by
\[
 \int_0^t\norm{\bar u}_{\mathcal B^J}
              \norm{\bar u}_{\mathcal B^{J+1}}\,\mathrm d\varrho
 \le t^{1/2}
     \norm{\bar u}_{\widetilde L^\infty_t\mathcal B^J}^{3/2}
     \norm{\bar u}_{L^1_t\mathcal B^{J+2}}^{1/2},
\]
and the Lorentz force is bounded by
$Ct\mathcal E_{B,J}(t)^2$. Applying \Cref{lem:tame-transport} to the
magnetic equation at index $J+1$ yields
\begin{equation}\label{eq:bg-high-B}
 \mathcal E_{B,J}(t)
 \le C_Je^{C_J\mathcal E_{u,J}(t)}
 \bigl[\mu\norm{b^{\mathrm{lo}}}_{\mathcal B^{J+1}}
       +\mathcal E_{B,J}(t)\mathcal E_{u,J}(t)\bigr].
\end{equation}

\hypertarget{not:reference-barriers}{}Take fixed constants
$A_B>4C_J\norm{b^{\mathrm{lo}}}_{\mathcal B^{J+1}}$ and
$A_u>4C_JA_B^2$. Bootstrap
$\mathcal E_{B,J}(t)<A_B\mu$ and
$\mathcal E_{u,J}(t)<A_u\mu^2t$.
The first inequality holds initially. Continuity and
\eqref{eq:bg-high-u} initialize the second: on a sufficiently short
interval, $C_Jt^{1/2}\mathcal E_{u,J}(t)\le1/2$, so that
$\mathcal E_{u,J}(t)\le2C_JA_B^2\mu^2t<A_u\mu^2t$.
Under both barriers, \eqref{eq:bg-high-u}--\eqref{eq:bg-high-B}
improve them to
\[
 \mathcal E_{u,J}(t)\le\tfrac12A_u\mu^2t,
 \qquad
 \mathcal E_{B,J}(t)\le\tfrac12A_B\mu
\]
when $\mu$ is sufficiently small, since $t\le1$. The first-contact
argument therefore gives
\begin{equation}\label{eq:bg-high-closed}
 \mathcal E_{u,J}(t)\le C_J\mu^2t,
 \qquad \mathcal E_{B,J}(t)\le C_J\mu
 \quad(0<t<\min\{1,T_{\max}\}).
\end{equation}
These bounds control the equal-order $H^J$ norms in the continuation
criterion of \Cref{lem:smooth-flow}, and hence $T_{\max}>1$.
Embeddings give \eqref{eq:bg-u} and \eqref{eq:bg-B-size}.

Finally, integrate the magnetic equation:
\[
 \bar B(t)-\mu b^{\mathrm{lo}}
 =\int_0^t\bigl(\bar B\cdot\nabla\bar u
                    -\bar u\cdot\nabla\bar B\bigr)\,\mathrm d\varrho.
\]
Because $J-1>\max\{r,\sigma\}$, the algebra estimate in
$\mathcal B^{J-1}$ bounds its
$\widetilde L^\infty_T\mathcal B^\sigma$ norm by
\[
 C_J\norm{\bar B}_{\widetilde L^\infty_T\mathcal B^{J+1}}
     \norm{\bar u}_{L^1_T\mathcal B^J}
 \le C_J\mu\mathcal E_{u,J}(T)
 \le C_\sigma\mu^3T.
\]
This proves \eqref{eq:bg-B-increment}.
\end{proof}

\subsection{Short-time perturbation bounds}

\begin{equation}\label{eq:perturbation-definition}
 v:=u-\bar u,\qquad h:=B-\bar B.
\end{equation}
Their initial data are
$(\lambda u_K^{\mathrm{hi}},0)$, and subtraction of the two exact
systems gives
\begin{align}
 \partial_tv-\Delta v
 &=-\mathbb P\bigl(\bar u\cdot\nabla v+v\cdot\nabla\bar u
                         +v\cdot\nabla v\bigr)\notag\\
 &\quad+\mathbb P\bigl(\bar B\cdot\nabla h+h\cdot\nabla\bar B
                         +h\cdot\nabla h\bigr),\notag\\
 \partial_th+(\bar u+v)\cdot\nabla h
 &=\bar B\cdot\nabla v+h\cdot\nabla v
       +h\cdot\nabla\bar u-v\cdot\nabla\bar B.
 \label{eq:perturbation-system}
\end{align}
All four fields are divergence-free.

\begin{proposition}[Endpoint and higher-order bounds]
\label[proposition]{prop:critical}
For every fixed integer $N>r+3$ there are
$\eta_N\in(0,1/4]$, $C>0$, and $C_N>0$ such that the following holds.
If the parameters in \eqref{eq:data-and-parameters} satisfy
$m+\mu+\tau\le\eta_N$ and $K\ge16$, then the classical solution
exists on $[0,T]$, where $T=\tau K^{-2}$, and
\begin{gather*}
 \norm{v}_{\widetilde L^\infty_T\mathcal B^{r-1}}\le Cm,
 \qquad
 \norm{v}_{L^1_T\mathcal B^{r+1}}\le Cm\tau,\\
 \norm{v}_{L^1_T\mathcal B^r}\le Cm\tau K^{-1},
 \qquad
 \norm{h}_{\widetilde L^\infty_T\mathcal B^r}\le C\mu m\tau.
\end{gather*}
Writing
\begin{equation}\label{eq:velocity-splitting}
 v_{\mathrm{lin}}:=e^{t\Delta}(\lambda u_K^{\mathrm{hi}}),
 \qquad v_{\mathrm{rem}}:=v-v_{\mathrm{lin}},
\end{equation}
one also has
\begin{equation}\label{eq:velocity-remainder-bound}
 \norm{v_{\mathrm{rem}}}_{\widetilde L^\infty_T\mathcal B^{r-1}}
 +\norm{v_{\mathrm{rem}}}_{L^1_T\mathcal B^{r+1}}
 \le Cm\tau(m+\mu+\tau),
\end{equation}
and the quantitative higher-order estimate
\begin{equation}\label{eq:velocity-high-quantitative}
 \norm{v}_{\widetilde L^\infty_T\mathcal B^{N-1}}
 +\norm{v}_{L^1_T\mathcal B^{N+1}}
 \le C_NmK^{N-r}.
\end{equation}
The constants $C$ are independent of $N$, while $\eta_N$ and $C_N$
may depend on $N$. All are independent of $K,\lambda,\mu,$ and $\tau$.
\end{proposition}

\begin{proof}
Fix $N$ before choosing the smallness threshold. For each fixed $K$,
let $[0,T_{\max})$ be the maximal $H^{N-1}\times H^{N-1}$ lifespan
from \Cref{lem:smooth-flow}. Persistence of all Sobolev orders
holds on its compact subintervals. Work on
$0<t<\min\{T,T_{\max}\}$. Define
\begin{equation}\label{eq:endpoint-running-norms}
\begin{aligned}
 X(t)&:=\norm{v}_{\widetilde L^\infty_t\mathcal B^{r-1}},
 &D(t)&:=\norm{v}_{L^1_t\mathcal B^{r+1}},\\
 Y(t)&:=\norm{h}_{\widetilde L^\infty_t\mathcal B^r},
 &E(t)&:=\norm{v_{\mathrm{rem}}}_{\widetilde L^\infty_t\mathcal B^{r-1}}
       +\norm{v_{\mathrm{rem}}}_{L^1_t\mathcal B^{r+1}}.
\end{aligned}
\end{equation}
These quantities are continuous on every compact interval inside the
smooth lifespan, by approximation with finite dyadic sums.

\smallskip
\noindent\emph{1. Endpoint inequalities.}
The Fourier support of $v_{\mathrm{lin}}$ lies where $|\xi|\simeq K$.
Thus, for a fixed $C_0\ge1$,
\begin{align}
 \norm{v_{\mathrm{lin}}}_{\widetilde L^\infty_t\mathcal B^{r-1}}
 &\le C_0m,\notag\\
 \norm{v_{\mathrm{lin}}}_{L^1_t\mathcal B^{r+1}}
 &\le C_0m(1-e^{-ctK^2})\le C_0m\tau.
 \label{eq:linear-truncated}
\end{align}
Consequently,
\begin{equation}\label{eq:linear-remainder-decomposition}
 X(t)\le C_0m+E(t),\qquad D(t)\le C_0m\tau+E(t).
\end{equation}
The heat estimate and the endpoint product laws applied to the first
equation in \eqref{eq:perturbation-system} give
\begin{equation}\label{eq:velocity-remainder-estimate}
 E(t)\le C\bigl[X(t)D(t)+\mu^2t\{X(t)+D(t)\}
                       +t(\mu+Y(t))Y(t)\bigr].
\end{equation}
Here the mixed velocity terms use
$\norm{\bar u}_{\widetilde L^\infty_t\mathcal B^{r-1}}
+\norm{\bar u}_{L^1_t\mathcal B^{r+1}}\le C\mu^2t$.
For each magnetic term, incompressibility places the derivative
outside a product of two $\mathcal B^r$ fields, so no higher norm of
$h$ occurs.

In the magnetic equation, retain transport by $\bar u+v$ on the left
and treat $h\cdot\nabla(\bar u+v)$ as stretching. The source terms
are bounded in $\mathcal B^r$ by
$C\mu\norm{v}_{\mathcal B^{r+1}}$ and
$C\mu\norm{v}_{\mathcal B^r}$, respectively. Hence
\Cref{lem:transport} yields, with $C_{\mathrm{tr}}\ge1$,
\begin{equation}\label{eq:magnetic-perturbation-estimate}
 Y(t)\le C_{\mathrm{tr}}
    e^{C_{\mathrm{tr}}(D(t)+\mu^2t)}
    \mu\left(D(t)+\norm{v}_{L^1_t\mathcal B^r}\right).
\end{equation}
The intermediate norm satisfies
\begin{equation}\label{eq:intermediate-velocity-bound}
 \norm{v}_{L^1_t\mathcal B^r}
 \le t^{1/2}X(t)^{1/2}D(t)^{1/2}.
\end{equation}

\smallskip
\noindent\emph{2. Closure of the endpoint bounds.}
Choose the fixed constants
\begin{equation}\label{eq:perturbation-barriers}
 A_v:=4C_0,\qquad A_h:=16C_{\mathrm{tr}}C_0,
\end{equation}
and bootstrap the three inequalities
\begin{equation}\label{eq:bootstrap-assumption}
 X(t)<A_vm,\qquad D(t)<A_vm\tau,\qquad
 Y(t)<A_h\mu m\tau.
\end{equation}
They hold on a nonempty initial interval because
$X(0)\le C_0m$ and $D(0)=Y(0)=E(0)=0$.
First substitute these barriers into
\eqref{eq:velocity-remainder-estimate}. Since $t\le\tau K^{-2}$ and
$m+\mu+\tau\le1$, this gives
\begin{equation}\label{eq:remainder-bootstrap}
 E(t)\le C_*m\tau(m+\mu+\tau),
\end{equation}
where $C_*$ depends only on the already fixed $A_v,A_h$ and the
endpoint constants. Choose the threshold small enough that
$C_*\eta_N\le C_0$ and the exponential in
\eqref{eq:magnetic-perturbation-estimate} is at most two under the
barriers. Equations \eqref{eq:linear-remainder-decomposition} and
\eqref{eq:remainder-bootstrap} then give the improved bounds
\[
 X(t)\le2C_0m=\tfrac12A_vm,
 \qquad D(t)\le2C_0m\tau=\tfrac12A_vm\tau.
\]
Use these improved bounds in \eqref{eq:intermediate-velocity-bound}
to obtain
$\norm{v}_{L^1_t\mathcal B^r}\le2C_0m\tau K^{-1}$.
Only now apply \eqref{eq:magnetic-perturbation-estimate}:
\begin{equation}\label{eq:magnetic-bootstrap}
 Y(t)\le4C_{\mathrm{tr}}C_0(1+K^{-1})\mu m\tau
       \le\tfrac12A_h\mu m\tau.
\end{equation}
Thus all three barriers improve strictly. The continuity argument
gives, for every $t<\min\{T,T_{\max}\}$,
\begin{gather}
 X(t)\le Cm,\quad D(t)\le Cm\tau,\quad
 \norm{v}_{L^1_t\mathcal B^r}\le Cm\tau K^{-1},\quad
 Y(t)\le C\mu m\tau,\label{eq:critical-truncated}\\
 E(t)\le Cm\tau(m+\mu+\tau).\label{eq:remainder-truncated}
\end{gather}
These estimates have used only endpoint norms on the truncated
interval; continuation to $T$ is established next.

\smallskip
\noindent\emph{3. Quantitative higher-order estimates.}
Define
\begin{equation}\label{eq:high-running-norms}
 \mathcal V_N(t)
 :=\norm{v}_{\widetilde L^\infty_t\mathcal B^{N-1}}
   +\norm{v}_{L^1_t\mathcal B^{N+1}},
 \qquad
 \mathcal H_N(t)
 :=\norm{h}_{\widetilde L^\infty_t\mathcal B^N}.
\end{equation}
The reference estimates give
\begin{equation}\label{eq:reference-high-size}
 \norm{\bar u}_{\widetilde L^\infty_T\mathcal B^{N-1}}
 +\norm{\bar u}_{L^1_T\mathcal B^{N+1}}\le C_N\mu^2T,
 \qquad
 \norm{\bar B}_{\widetilde L^\infty_T\mathcal B^{N+1}}
 \le C_N\mu.
\end{equation}
We decrease $\eta_N$ so that \Cref{lem:background} applies at all
indices in this display. We also impose
\begin{equation}\label{eq:high-absorption-smallness}
 C_N\norm{u}_{\widetilde L^\infty_t\mathcal B^{r-1}}
 \le C_N(m+\mu^2T)\le\tfrac14,
\end{equation}
where $C_N$ includes the constant in the tame parabolic estimate.
This restriction involves no higher norm of the perturbation.

Write the velocity equation as
\[
 \partial_tv-\Delta v+\mathbb P(u\cdot\nabla v)
 =-\mathbb P(v\cdot\nabla\bar u)
  +\mathbb P(\bar B\cdot\nabla h+h\cdot\nabla\bar B
                                          +h\cdot\nabla h).
\]
The estimate in \Cref{lem:tame-parabolic} absorbs the term
$C_N\norm{u}_{\widetilde L^\infty_t\mathcal B^{r-1}}
\norm{v}_{L^1_t\mathcal B^{N+1}}$ by
\eqref{eq:high-absorption-smallness}. The remaining convection terms
are bounded by
\begin{align*}
 &C_N\int_0^t
   \bigl(\norm{v}_{\mathcal B^{r+1}}
         +\norm{\bar u}_{\mathcal B^{r+1}}\bigr)
   \mathcal V_N(\varrho)\,\mathrm d\varrho\\
 &\qquad+C_Nm\tau
       \norm{\bar u}_{\widetilde L^\infty_T\mathcal B^{N-1}}
       +C_Nm\norm{\bar u}_{L^1_T\mathcal B^{N+1}}.
\end{align*}
The divergence-form tame product bound gives
\begin{align}
 &\norm{\bar B\cdot\nabla h+h\cdot\nabla\bar B+h\cdot\nabla h}
             _{L^1_t\mathcal B^{N-1}}\notag\\
 &\quad\le C_NT\mu\mathcal H_N(t)
       +C_NT\norm{\bar B}_{\widetilde L^\infty_T\mathcal B^N}Y(t)
 \le C_NT\mu\mathcal H_N(t)+C_NT\mu^2m\tau.
 \label{eq:high-v-magnetic}
\end{align}
Since
$\norm{v(0)}_{\mathcal B^{N-1}}\le C_NmK^{N-r}$, all terms
depending only on the reference flow and the endpoint bounds can be
included in this initial size. We obtain
\begin{equation}\label{eq:high-v-ledger-closed}
 \begin{split}
 \mathcal V_N(t)
 &\le C_NmK^{N-r}+C_NT\mu\mathcal H_N(t)\\
 &\quad+C_N\int_0^t
       \bigl(\norm{v}_{\mathcal B^{r+1}}
             +\norm{\bar u}_{\mathcal B^{r+1}}\bigr)
       \mathcal V_N(\varrho)\,\mathrm d\varrho.
 \end{split}
\end{equation}

For the magnetic equation, \Cref{lem:tame-transport} at index $N$
gives
\begin{align}
 \mathcal H_N(t)
 &\le C_N\exp\!\left(C_N\norm{u}_{L^1_t\mathcal B^{r+1}}\right)
 \biggl[\norm{\bar B\cdot\nabla v-v\cdot\nabla\bar B}
                         _{L^1_t\mathcal B^N}\notag\\
 &\hspace{11em}
  +Y(t)\norm{u}_{L^1_t\mathcal B^{N+1}}\biggr].
 \label{eq:high-h-raw}
\end{align}
The exponential is bounded, since
$\norm{u}_{L^1_t\mathcal B^{r+1}}\le Cm\tau+C_N\mu^2T$.
The high-order factors in the source and stretching terms obey
\begin{align}
 \norm{\bar B\cdot\nabla v}_{L^1_t\mathcal B^N}
 &\le C_N\mu\mathcal V_N(t)+C_N\mu m\tau,
 \label{eq:high-h-source-one}\\
 \norm{v\cdot\nabla\bar B}_{L^1_t\mathcal B^N}
 &\le C_N\mu\norm{v}_{L^1_t\mathcal B^N}
       +C_N\mu\norm{v}_{L^1_t\mathcal B^r}\notag\\
 &\le C_N\mu\mathcal V_N(t)+C_N\mu m\tau K^{-1},
 \label{eq:high-h-source-two}\\
 Y(t)\norm{u}_{L^1_t\mathcal B^{N+1}}
 &\le C\mu m\tau\mathcal V_N(t)+C_N\mu^3m\tau T.
 \label{eq:high-h-stretch}
\end{align}
The first line assigns the high derivative either to $v$ or to the
smooth reference field. In the second line, the two terms arise from
$\norm{\bar B}_{\mathcal B^{r+1}}
\norm{v}_{\mathcal B^N}$ and
$\norm{\bar B}_{\mathcal B^{N+1}}\norm{v}_{\mathcal B^r}$.
Thus no $\mathcal B^{N+1}$ norm of $h$ is required. Combining these
bounds with \eqref{eq:high-h-raw} yields
\begin{equation}\label{eq:high-h-ledger-closed}
 \mathcal H_N(t)\le C_N\mu m\tau+C_N\mu\mathcal V_N(t).
\end{equation}

Substitution in \eqref{eq:high-v-ledger-closed} produces
$C_NT\mu^2\mathcal V_N(t)$, which is absorbed by decreasing
$\eta_N$ once more. The remaining integral coefficient has total
size at most $Cm\tau+C_N\mu^2T$. Gronwall's inequality gives the
quantitative bounds
\begin{equation}\label{eq:high-persistence-closed}
 \mathcal V_N(t)\le C_NmK^{N-r},
 \qquad
 \mathcal H_N(t)\le C_N\mu mK^{N-r}
 \quad(0<t<\min\{T,T_{\max}\}).
\end{equation}
Here $C_N$ is independent of $K$; the explicitly displayed power of
$K$ accounts for the high-frequency initial norm.

\smallskip
\noindent\emph{4. Continuation.}
For each fixed $K$, \eqref{eq:high-persistence-closed} and
\eqref{eq:reference-high-size} uniformly control $u$ and $B$ in
$H^{N-1}$ up to $\min\{T,T_{\max}\}$. Since $N-1>r+1$,
\Cref{lem:smooth-flow} excludes $T_{\max}\le T$. Letting $t\uparrow T$
in \eqref{eq:critical-truncated}, \eqref{eq:remainder-truncated}, and
\eqref{eq:high-persistence-closed} proves the proposition. The datum
belongs to every Sobolev space, so persistence of regularity makes
the solution classical throughout $[0,T]$.
\end{proof}

The higher-order bound also controls the velocity in the topology of
the theorem. Suppose $0<\alpha<s-r$, put $q=s-1-\alpha$, and fix
$N>\max\{r+3,q+1\}$ before choosing the parameters. Then
$r-1<q<N-1$, and \eqref{eq:CL-interpolation} with
$\gamma=(q-r+1)/(N-r)$ gives
\[
 \norm{v}_{\widetilde L^\infty_T\mathcal B^q}
 \le C_N
   \norm{v}_{\widetilde L^\infty_T\mathcal B^{r-1}}^{1-\gamma}
   \norm{v}_{\widetilde L^\infty_T\mathcal B^{N-1}}^\gamma
 \le C_NmK^{q-r+1}=C_N\lambda.
\]
Together with \Cref{lem:background} and
$\mathcal B^q\hookrightarrow H^q$, this proves
\begin{equation}\label{eq:target-velocity-small}
 \sup_{0\le t\le T}\norm{u(t)}_{H^q}
 \le C_N(\lambda+\mu^2T).
\end{equation}

\subsection{Localized expansion of the magnetic field}

We compare the exact perturbation $h$ with $\mathcal M_K$ in
$\mathcal B^r$, then recover an $H^s$ lower bound from the selected
frequency packet.

\begin{proposition}[Expansion at the selected frequency]
\label[proposition]{prop:expansion}
Fix an integer $N>r+3$ and $L>0$. There exists
$\eta_{N,L}\in(0,\eta_N]$ such that, if the parameters in
\eqref{eq:data-and-parameters} satisfy
$m+\mu+\tau\le\eta_{N,L}$, then
\begin{equation}\label{eq:expansion-error}
 \begin{split}
 K^r\norm{\mathcal Q_K\bigl(B(T)-\mathcal M_K(T)\bigr)}_2
 &\le C\mu m\tau(m+\mu+\tau+K^{-1})\\
 &\quad+C_L\mu^3TK^{-L}.
 \end{split}
\end{equation}
The constants are independent of $K,\lambda,\mu,\tau$.
\end{proposition}

\begin{proof}
Choose $\eta_{N,L}\le\eta_N$ small enough that
\Cref{lem:background} also applies at index $r+L$.
Set
\begin{equation}\label{eq:magnetic-remainder-definition}
 \mathcal R_K:=h-\mathcal M_K.
\end{equation}
Subtracting
$\partial_t\mathcal M_K=\mu b^{\mathrm{lo}}\cdot\nabla v_{\mathrm{lin}}$
from \eqref{eq:perturbation-system} gives the exact equation
\begin{equation}\label{eq:nonlinear-remainder-equation}
 \begin{split}
 \partial_t\mathcal R_K+(\bar u+v)\cdot\nabla\mathcal R_K
 &=\mu b^{\mathrm{lo}}\cdot\nabla v_{\mathrm{rem}}
   +(\bar B-\mu b^{\mathrm{lo}})\cdot\nabla v\\
 &\quad+h\cdot\nabla v+h\cdot\nabla\bar u
        -v\cdot\nabla\bar B\\
 &\quad-(\bar u+v)\cdot\nabla\mathcal M_K,
 \qquad \mathcal R_K(0)=0.
 \end{split}
\end{equation}
We use the already proved bound for $h$ in these source terms.
Transport of $\mathcal R_K$ remains on the left, so no derivative
of $\mathcal R_K$ is estimated as a source.

The pure transport estimate in \Cref{lem:transport} bounds
$\norm{\mathcal R_K}_{\widetilde L^\infty_T\mathcal B^r}$
by the sum of the $L^1_T\mathcal B^r$ norms of the sources.
Its exponential factor is bounded because
$\norm{\bar u+v}_{L^1_T\mathcal B^{r+1}}
\le C(\mu^2T+m\tau)$.
Applying the algebra property of $\mathcal B^r$,
\Cref{lem:background,lem:M-lower,prop:critical} gives the seven
bounds recorded in \Cref{tab:remainder-sources}.

\begin{table}[htbp]
\centering
\captionsetup{font=footnotesize,labelfont=bf}
\small
\renewcommand{\arraystretch}{1.16}
\begin{tabularx}{\textwidth}{@{}Y L{0.31\textwidth}@{}}
\toprule
Source $g$ in \eqref{eq:nonlinear-remainder-equation}
& Bound for $\norm{g}_{L^1_T\mathcal B^r}/(\mu m\tau)$\\
\midrule
$\mu b^{\mathrm{lo}}\cdot\nabla v_{\mathrm{rem}}$
& $C(m+\mu+\tau)$\\
$(\bar B-\mu b^{\mathrm{lo}})\cdot\nabla v$
& $C\mu^2T$\\
$h\cdot\nabla v$
& $Cm\tau$\\
$h\cdot\nabla\bar u$
& $C\mu^2T$\\
$-v\cdot\nabla\bar B$
& $CK^{-1}$\\
$-v\cdot\nabla\mathcal M_K$
& $Cm\tau$\\
$-\bar u\cdot\nabla\mathcal M_K$
& $C\mu^2\tau K^{-1}$\\
\bottomrule
\end{tabularx}
\caption{Integrated source bounds relative to the leading magnetic
scale. Each row estimates the entire source shown in the first column.}
\label{tab:remainder-sources}
\end{table}

In particular, the derivative on the smooth reference magnetic field
costs no factor of $K$:
\[
 \norm{v\cdot\nabla\bar B}_{L^1_T\mathcal B^r}
 \le C\mu\norm{v}_{L^1_T\mathcal B^r}
 \le C\mu m\tau K^{-1}.
\]
For transport of the explicit response, its one higher norm is
compensated by the same time-interpolation gain:
\begin{align*}
 &\norm{(\bar u+v)\cdot\nabla\mathcal M_K}_{L^1_T\mathcal B^r}\\
 &\qquad\le C\bigl(\mu^2T+m\tau K^{-1}\bigr)
             K\mu m\tau
 =C\mu m\tau\bigl(\mu^2\tau K^{-1}+m\tau\bigr).
\end{align*}
Since $m+\mu+\tau\le1$ and $T\le1$, summing the table yields
\begin{equation}\label{eq:nonlinear-remainder-bound}
 \norm{\mathcal R_K}_{\widetilde L^\infty_T\mathcal B^r}
 \le C\mu m\tau(m+\mu+\tau+K^{-1}).
\end{equation}

Finally, $B=\bar B+\mathcal M_K+\mathcal R_K$, while
$\mathcal Q_K(\mu b^{\mathrm{lo}})=0$. The smooth reference tail
satisfies
\begin{equation}\label{eq:bg-tail}
 \begin{split}
 K^r\norm{\mathcal Q_K\bar B(T)}_2
 &\le C_LK^{-L}
       \norm{\bar B(T)-\mu b^{\mathrm{lo}}}_{\mathcal B^{r+L}}\\
 &\le C_L\mu^3TK^{-L}.
 \end{split}
\end{equation}
Combine this with \eqref{eq:nonlinear-remainder-bound} and
\eqref{eq:loc-Bs} to obtain \eqref{eq:expansion-error}.
\end{proof}

\section{Proof of norm inflation}\label{sec:norm-inflation}

\begin{proof}[Proof of \Cref{thm:main}]
It suffices to prove the result for $r-1<q<s-1$, since both
initial and evolved velocity norms decrease when the Sobolev index
is lowered. We therefore assume
\[
 0<\alpha:=s-1-q<s-r.
\]
Fix an integer $N>\max\{r+3,q+1\}$ and a tail order
$L>s-r+2$. These indices are chosen before $K$.
Choose
\begin{equation}\label{eq:delta-choice}
 0<\delta<\frac{\alpha}{4},
\end{equation}
and, for dyadic $K\to\infty$, set
\begin{equation}\label{eq:parameter-choice}
 \lambda=\mu=\tau=K^{-\delta},\qquad
 T=K^{-2-\delta},\qquad
 m=K^{-(s-r-\alpha+\delta)}.
\end{equation}
Since $s-r-\alpha>0$,
$m+\mu+\tau\to0$. Thus \Cref{prop:critical,prop:expansion} apply
for all sufficiently large $K$, and the solution exists beyond $T$.
The initial norms satisfy
\begin{equation}\label{eq:small-initial-data}
 \norm{u_0}_{H^q}+\norm{B_0}_{H^s}
 \le C(\lambda+\mu)\longrightarrow0.
\end{equation}
Moreover, \eqref{eq:target-velocity-small} gives
\begin{equation}\label{eq:small-evolved-velocity}
 \sup_{0\le t\le T}\norm{u(t)}_{H^q}
 \le C_N\bigl(K^{-\delta}+K^{-2-3\delta}\bigr)
 \longrightarrow0.
\end{equation}

It remains to verify that the error is negligible relative to the
positive response. The first relative error in
\eqref{eq:expansion-error} tends to zero. For the reference tail,
\[
 \frac{\mu^3TK^{-L}}{\mu m\tau}
 =K^{s-r-\alpha-2-L-\delta}\longrightarrow0
\]
by the choice of $L$. Hence
\Cref{lem:M-lower,prop:expansion} and the triangle inequality imply
\[
 K^r\norm{\mathcal Q_KB(T)}_2
 \ge\frac{c_0}{2}\mu m\tau
\]
for all sufficiently large $K$. By \eqref{eq:loc-Hs},
\begin{equation}\label{eq:norm-inflation-final}
 \norm{B(T)}_{H^s}
 \ge cK^{s-r}\mu m\tau
 =cK^{\alpha-3\delta}\longrightarrow\infty.
\end{equation}
The exponent is positive by \eqref{eq:delta-choice}.
Given $\varepsilon>0$ and $M>0$, choose a sufficiently large dyadic
$K$ so that \eqref{eq:small-initial-data} and
\eqref{eq:small-evolved-velocity} are below $\varepsilon$,
\eqref{eq:norm-inflation-final} exceeds $M$, and $T<\varepsilon$.
The profiles are real, divergence-free, and Schwartz, completing
the proof.
\end{proof}

\begin{proof}[Proof of \Cref{cor:no-continuous-solution-map}]
\hypertarget{not:corollary-sequence}{}Suppose such a map $\Phi$ exists. Agreement with the zero classical
solution implies $\Phi(0)=0$. Apply \Cref{thm:main} with
$\varepsilon=\min\{T_0/2,1/n\}$ and $M=n$, $n\ge1$. It gives
Schwartz data $z_{0,n}\to0$ in $\mathcal X_{q,s}$, magnetic solutions
$B_n$, and times $0<T_n<\min\{T_0/2,1/n,T_{\max}(z_{0,n})\}$
such that $\norm{B_n(T_n)}_{H^s}>n$. For large $n$, $z_{0,n}\in\mathcal U$,
and agreement during the classical lifespan gives
\[
 \norm{\Phi(z_{0,n})}_{C([0,T_0];\mathcal X_{q,s})}
 \ge\norm{B_n(T_n)}_{H^s}>n.
\]
This contradicts continuity at zero.
\end{proof}

\section*{Acknowledgments}
The work of Y. H. Wang was partially supported by the National Natural
Science Foundation of China grant 12401274.
\section*{Data availability statement}
No new data were created or analyzed in this study.
\section*{Conflict of interest}
The author declares no conflict of interest.

\section*{Declaration of generative AI use}
The research question was formulated by the author on the basis of
the conclusions in
\cite{FeffermanEtAl2014,FeffermanEtAl2017,Li2020}.
By examining the structure of the equations, the author conjectured
the principal magnetic norm-inflation conclusion now stated in
\Cref{thm:main}.

The initial draft of the proof was generated by ChatGPT 5.6.
The manuscript text, including the mathematical arguments and their
exposition, was generated entirely by ChatGPT-family models in
response to simple prompts supplied by the author. The author directed
the organization through prompts formulated from a reader's
perspective, and Codex carried out the corresponding revisions.

The first version that the author considered reasonably readable
comprised 30 pages and was produced on 13 August 2026. The author
nevertheless considered it too long. More recently, the author used
ChatGPT 6 Pro to streamline the arguments and shorten the
presentation, resulting in the present version.

The author reviewed the mathematical correctness of the arguments
and the manner in which they were presented, and judged that the
resulting manuscript required only limited further manual adjustment.
The author takes responsibility for the mathematical claims, proofs,
and final text.

The author considers this process worth documenting. In the author's view, this experience suggests that, with appropriate guidance, current ChatGPT models can organize a PDE research article of about twenty pages into an exposition accessible to readers with a basic background in modern PDE theory.

\end{document}